\documentclass[5p,preprint,10pt]{elsarticle}
\usepackage{lineno}
\usepackage{mathtools}
\usepackage{algorithm2e}
\usepackage{float}
\usepackage{fbox}
\usepackage{tikz}
\usetikzlibrary{shapes,arrows,positioning}

\newcommand{\dd}{\, \mathrm d}
\newcommand{\R}{\mathbb R}

\newcommand{\cG}{\mathcal G}

\newcommand{\sC}{\boldsymbol{\mathsf C}}
\newcommand{\brho}{\boldsymbol\rho}

\newcommand{\x}{\boldsymbol x}
\newcommand{\y}{\boldsymbol y}

\newcommand{\bdelta}{{\boldsymbol \delta}}
\newcommand{\cc}{{\boldsymbol c}}

\newcommand{\hxxiC}{\boldsymbol{\hat \xi^{\mathsf C}}}
\newcommand{\hxxiCo}{\boldsymbol{\hat \xi^{\mathsf C}}_{\!\!\!\!\!1}\,}
\newcommand{\hxxiCt}{\boldsymbol{\hat \xi^{\mathsf C}}_{\!\!\!\!\!2}\,}
\newcommand{\xxxic}{\boldsymbol{\xi^{\mathsf c}}}
\newcommand{\xxxico}{\boldsymbol{\xi^{\mathsf c}}_{\!\!\!1}}
\newcommand{\xxxict}{\boldsymbol{\xi^{\mathsf c}}_{\!\!\!2}}
\newcommand{\txxxic}{\boldsymbol{\tilde \xi^{\mathsf c}}}
\newcommand{\txxxico}{\boldsymbol{\tilde \xi^{\mathsf c}}_{\negthickspace\negthickspace 1}}
\newcommand{\txxxict}{\boldsymbol{\tilde \xi^{\mathsf c}}_{\negthickspace\negthickspace 2}}
\newcommand{\mX}{{\boldsymbol{\mathcal X}}}
\newcommand{\ccG}{\boldsymbol \cG}

\newcommand{\oomega}{\boldsymbol \omega}
\newcommand{\m}{\boldsymbol m}
\newcommand{\V}{\boldsymbol V\!}
\newcommand{\GG}{\boldsymbol G}
\newcommand{\ct}{\boldsymbol t}
\newcommand{\cv}{\boldsymbol v}

\usepackage[utf8]{inputenc}
\usepackage[]{algorithm2e}
\usepackage{comment}
\usepackage{graphicx}
\usepackage{amssymb,amsthm,mathtools,amsmath}
\usepackage{hyperref}
\usepackage[nameinlink]{cleveref}
\usepackage{todonotes}
\usepackage{dsfont}
\usepackage{algorithm2e}
\usepackage{capt-of}
\usepackage{mathtools, cuted}
\usepackage{lipsum, color}
\usepackage{nicefrac}
\usepackage{algorithmic}
\usepackage{letltxmacro}
\LetLtxMacro{\oldalgorithmic}{\algorithmic}
\LetLtxMacro{\endoldalgorithmic}{\endalgorithmic}
\renewenvironment{algorithmic}[1][0]{%
  \hrulefill\par
  \oldalgorithmic[#1]}
  {\endoldalgorithmic\par
   \vspace*{-.5\baselineskip}
   \hrulefill\par
  }
\usepackage{placeins}
\newenvironment{talign*}
 {\csname align*\endcsname}
 {\endalign}

\theoremstyle{plain}
\newtheorem{theorem}{Theorem}[section]

\newtheorem{defi}[theorem]{Definition}

\newtheoremstyle{boldremark}
    {\dimexpr\topsep/2\relax} 
    {\dimexpr\topsep/2\relax} 
    {}          
    {}          
    {\bfseries\itshape} 
    {.}         
    {.5em}      
    {}          
\theoremstyle{boldremark}    
 
\journal{-}

\begin{document}

\begin{frontmatter}

\title{Exact Method of Moments for multi-dimensional population balance equations}
\cortext[cor1]{Corresponding author}

\author[label0,label1]{Adeel Muneer}
\author[label3]{Tobias Schikarski}
\author[label0,label1]{Lukas Pflug\corref{cor1}}
\ead{lukas.pflug@fau.de}

\address[label0]{Friedrich-Alexander-Universität Erlangen-Nürnberg (FAU), Competence Unit for Scientific Computing (CSC), Martensstr. 5a, 91058 Erlangen, Germany}
\address[label1]{Friedrich-Alexander-Universität Erlangen-Nürnberg (FAU), Chair of Applied Mathematics (Continuous Optimization), Department of Mathematics, Cauerstraße 11, 91058 Erlangen, Germany}
\address[label3]{Friedrich-Alexander-Universität Erlangen-Nürnberg (FAU), Institute of Particle Technology (LFG), Cauerstrasse 4, 91058, Erlangen, Germany}

\begin{abstract}
The unique properties of anisotropic and composite particles are increasingly being leveraged in modern particulate products. However, tailored synthesis of particles characterized by multi-dimensional dispersed properties remains in its infancy and few mathematical models for their synthesis exist. Here, we present a novel, accurate and highly efficient numerical approach to solve a multi-dimensional population balance equation, based on the idea of the exact method of moments for nucleation and growth \cite{pflug2020emom}. The transformation of the multi-dimensional population balance equation into a set of one-dimensional integro-differential equations allows us to exploit accurate and extremely efficient numerical schemes that markedly outperform classical methods (such as finite volume type methods) which is outlined by convergence tests. Our approach not only provides information about complete particle size distribution over time, but also offers insights into particle structure. The presented scheme and its performance is exmplified based on coprecipitation of nanoparticles. For this process, a generic growth law is derived and parameter studies as well as convergence series are performed.
\end{abstract}

\begin{keyword}
multi-dimensional \sep exact method of moments \sep population balance equation \sep fixed point equation \sep integro-differential equation \sep  growth kinetics \sep synthesis \sep nanoparticle \sep nonlocal conservation laws \sep method of characteristics
\end{keyword}
\end{frontmatter}

\section{Introduction and problem definition}
The modeling and efficient numerical approximations of multi-dimensional (MD) nanoparticle (NP) synthesis are increasingly required as the properties of anisotropic particles play a major role in various applications, including bio-nanosensors or catalysts \cite{pmid17722542,C2CP43642F, PEUKERT20151}. In recent years, automated analysis of MD particle shape distributions -- essential to validate and calibrate MD process models --  has been demonstrated for a range of composite and anisotropic NPs \cite{wawra2018determination,meincke2022determination,lopez2022multidimensional}, unlocking the possibility of using predictive modeling for MD-NP synthesis. 

To this end, we here extend the recently derived \textit{exact Method of Moments} (eMoM\cite{pflug2020emom}) to MD population balance equations (PBEs). As an example, we study a class of balance laws describing coprecipitation of particles characterized by their size and composition, such as the seeded growth of nanoalloys. However, the scheme can be applied more generally and is also applicable to e.g., the modeling of shape anisotropic growth (see e.g., the kinetics analyzed in \cite{braatz2002particle,SCHIELE202354}). In general, we study the following MD-PBE:

\begin{defi}[Multi-dimensional population balance equation] 
The evolution of the particle population $q$ can be macroscopically described by the following MD-PBE: 
\begin{align}
\begin{aligned}
q_t(t,\x)+\nabla_{\x} (\ccG(\cc(t),\x)q(t,\x))&=0,\\ q(0,\x)&=q_0(\x), 
\end{aligned}\label{eqn:model}
\end{align}
 with $\x \in \mathcal X \subset \R^n$, $t\in [0,T]$. The concentration of the $i$-th educt species, i.e. $\cc_i$, is defined by:
\begin{align}
\cc_i(t)&=\tfrac{1}{\mathcal V}\m_i(t)-\tfrac{\brho_i}{\mathcal V}\iint\limits_{\mX}\V_i(\x) q(t,\x) \, \mathrm d \x \label{eqn:COM1},
\end{align}
where $q$ denotes the particle number density, $t$ the process time, $\x$ the shape parametrization of the dispersed phase, $\ccG$ the concentration-dependent growth rate in the first argument and shape in the second argument, $q_0$ the initial particle number density, $\cc_i$ the concentration of the $i$-th educt species, $\m_i$ the total mass of the $i$-th precipitated component, $\mathcal V$ the reactor volume, $\brho_i$ the density of the $i$-th species in the NPs, $\mX$ the set of admissible particle shapes, and $\V_i(\x)$ the volume of the $i$-th component in a particle of shape $\x$.   
\end{defi}

Conservation of mass (\cref{eqn:COM1}) for each component couples the PBE solution $q$ and the concentrations $c_i$, rendering the PBE \cref{eqn:model} and \eqref{eqn:COM1} a (multi-dimensional) nonlocal conservation law \cite{keimer2018existence}. In this manuscript, we reformulate the MD-PBE purely in terms of the concentrations, i.e., the driving forces. This idea was recently applied to the one-dimensional case in \cite{pflug2020emom}. The advantage of such a reformulation lies in reducing the MD-PBE \cref{eqn:model} to a set of one-dimensional integro-differential equations \cref{eqn:COM2}. The resulting equations prescribing the time-evolution of the concentrations can then be efficiently discretized and numerically approximated. For the most simple case involving two-chemical components and the evolution of composite particles, this reformulation reduces the need for numerically approximating one three-dimensional function, here the number density function $q$, to approximating two one-dimensional functions, i.e., the evolution of the two concentrations $c_i$.

\section{Multi-dimensional exact method of moments}\label{sec:derivation_eMoM}
Our aim is to derive a solely concentration-dependent ($\cc$) equation based on \cref{eqn:COM1}. Using the method of characteristics, see e.g. \cite{keimer2018existence}, for a given concentration $\cc$, the solution of \cref{eqn:model} can be stated as follows:
\begin{align}
q(t,\x)&=q_0(\xxxic [t,\x](0))\text{det}(D_2\xxxic [t,\x](0)) \label{eqn:SOL1},
\end{align}
with the characteristics satisfying:
\begin{align}
\begin{aligned}
\partial_3 \xxxic [t,\x](\tau)&=\ccG(\cc(\tau),\xxxic [t,\x](\tau)),\\
\xxxic [t,\x](t)&=\x, 
\end{aligned}\label{eqn:CHAR1}
\end{align}
for every $(\tau,t,\x)\in [0,T]^2 \times \mX$ for $T \in \R_{>0}$ and with $\mX \subset \R^2$.
As the characteristics depend on the concentration, we have added the concentration $\cc$ as a superscript to illustrate this dependency. By plugging the solution formula \cref{eqn:SOL1} into \cref{eqn:COM1}, and integrating by substitution using $\x=\xxxic[0,\y](t)$, we end up -- assuming $\V_i(\x) = 0 \, \forall \x \notin \mX$ -- with the following integro-differential equation:
\begin{defi}[Integral fixed-point-problem for the concentrations $\cc$]
\begin{align}
\cc_i(t)&=\tfrac{1}{\mathcal V} \m_i(t)-\tfrac{\brho_i}{\mathcal V}\iint_{\mX} \V_i(\xxxic[0,\y](t))  q_0\left(\y\right)\, \mathrm d \y \label{eqn:COM2}.
\end{align}
\end{defi}

The idea of eMoM (\cite{pflug2020emom}) is now to numerically approximate \cref{eqn:COM2} instead of \cref{eqn:model} and \eqref{eqn:COM1}. With $\cc$ as a solution of \cref{eqn:COM2}, we can \textit{post hoc} evaluate the PBE solution $q$ using \cref{eqn:SOL1}. Thus, to obtain the full PBE solution $q$, we only need to compute a small number (the number of components involved) of time-dependent scalar quantities, i.e., the concentration $\cc$. Thus, the reformulation \cref{eqn:COM2} of \cref{eqn:model} and \eqref{eqn:COM1} is highly advantageous whenever \cref{eqn:CHAR1} can be solved analytically for an arbitrary but given concentration $\cc$, which is indeed possible for a distinct class of growth law functions. In the following, we first introduce rather general growth kinetics for coprecipitation processes, and then derive the analytical solution of the corresponding characteristics equation.

\section{Kinetics of coprecipitation}
In this section, we derive -- based on generic physical assumptions -- the growth kinetics of NPs in a coprecipitation process, in which two educt species simultaneously assemble in one particle ensemble. This approach is, e.g., in line with the current process understanding of synthesized nanoalloys \cite{rioux2015seeded}, but also applies to the synthesis of multi-component battery cathode materials such as nickel-manganese-cobalt hydroxide \cite{coprep_nmcoh}.
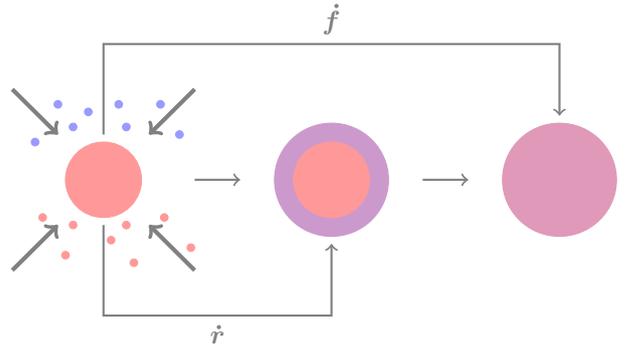
\begin{figure}
    \centering
\begin{tikzpicture}
\filldraw[color=blue!40, fill=blue!40!white](0.75,1.0) circle (0.05); 
\filldraw[color=blue!40, fill=blue!40!white](1.0,0.6) circle (0.05); 
\filldraw[color=blue!40, fill=blue!40!white](0.2,1.0) circle (0.05); 
\filldraw[color=blue!40, fill=blue!40!white](0.3,0.7) circle (0.05); 
\filldraw[color=blue!40, fill=blue!40!white](-0.2,0.9) circle (0.05); 
\filldraw[color=blue!40, fill=blue!40!white](-0.4,0.7) circle (0.05); 
\filldraw[color=blue!40, fill=blue!40!white](-0.9,0.5) circle (0.05); 
\filldraw[color=blue!40, fill=blue!40!white](-0.6,1.0) circle (0.05); 
\filldraw[color=red!40, fill=red!40!white](1.15,-0.9) circle (0.05); 
\filldraw[color=red!40, fill=red!40!white](0.8,-0.5) circle (0.05); 
\filldraw[color=red!40, fill=red!40!white](0.4,-1.1) circle (0.05); 
\filldraw[color=red!40, fill=red!40!white](0.3,-0.6) circle (0.05); 
\filldraw[color=red!40, fill=red!40!white](0.1,-0.8) circle (0.05); 
\filldraw[color=red!40, fill=red!40!white](-0.4,-0.6) circle (0.05); 
\filldraw[color=red!40, fill=red!40!white](-0.5,-1.0) circle (0.05); 
\filldraw[color=red!40, fill=red!40!white](-0.8,-0.5) circle (0.05); 
\filldraw[color=red!40, fill=red!40!white](0,0) circle (0.5);
\filldraw[color=red!50!blue!40!white, fill=red!50!blue!40!white](3,0) circle (0.75);
\filldraw[color=red!40, fill=red!40!white](3,0) circle (0.5);
\filldraw[color=red!70!blue!40!white, fill=red!70!blue!40!white](6,0) circle (0.75);
\draw[gray, ultra thick,->] (1.2,1.2) -- (0.6,0.6); 
\draw[gray, ultra thick,->] (-1.2,1.2) -- (-0.6,0.6);
\draw[gray, ultra thick,->] (-1.2,-1.2) -- (-0.6,-0.6); 
\draw[gray, ultra thick,->] (1.2,-1.2) -- (0.6,-0.6);
\draw[gray, thick,->] (0,-0.6) -- (0,-1.8) -- node[midway,below]{$\boldsymbol{\dot{r}}$}(3,-1.8) -- (3,-0.85);
\draw[gray, thick,->] (0,0.6) -- (0,1.8) -- node[midway,above]{$\boldsymbol{\dot{f}}$}(6,1.8)-- (6,0.85);
\draw[gray, thick,->] (1.2,0) -- (1.8,0); 
\draw[gray, thick,->] (4.2,0)  --  (4.8,0); 
\end{tikzpicture}
\caption{Sketch of the change in size ($\dot r$) and composition $(\dot f)$ in the growth of alloy nanoparticles}
\label{fig:growth_kinetics_sketch}
\end{figure}

We assume that the change in NP radius is given by the sum of two growth rates, which depend on the radius to the exponent $n$ and the two concentrations $\cc_1$ and $\cc_2$, i.e.:
\begin{align}
\dot{r}(t)= \left(\GG_1(\cc_1(t))+\GG_2(\cc_2(t))\right) r(t)^n, \label{eqn:GROW1}
\end{align}
 where $\GG_1$ and $\GG_2$ denote the growth kinetics of the two species. For the change in NP size, we can thus model, e.g. size-independent growth, i.e. $n=0$ \cite{Thanh2014}, or diffusion-limited growth, i.e., $n=-1$ \cite{schikarski2021cej,schikarski2022predictive}.\\
For a particle with a given radius ($r>0$), we further define the volume fraction or composition $f$ and its derivative as: 
\begin{align}
\begin{aligned}
f(t)&=\tfrac{\cv_1(t)}{\cv_1(t)+\cv_2(t)}, \\ \dot{f}(t)&=\tfrac{\dot{\cv_1}(t)}{\cv_1(t)+\cv_2(t)}-f(t) \tfrac{\dot{\cv_1}(t)+\dot{\cv_2}(t)}{\cv_1(t)+\cv_2(t)}.
\end{aligned}\label{eqn:GROW3}
\end{align}
Here, $\cv_1(t)$ and $\cv_2(t)$ are the volumes of the first and second chemical components in the NP at time $t$, respectively. Clearly, $\cv_1(t)+\cv_2(t) = \frac{4}{3}\pi r(t)^3$, and thus $\dot \cv_1(t)+\dot \cv_2(t)=4\pi r^2 \dot r = 4\pi r^{n+2}\left(\GG_1(\cc_1(t))+\GG_2(\cc_2(t))\right)$. As $\cv_1(t)$ and $\cv_2(t)$ can only change due to $\GG_1(\cc_1(t))$ and $\GG_2(\cc_2(t))$, respectively, we can assign
$\dot \cv_1(t) =4 \pi r^{n+2} \GG_1(\cc_1(t))$, and $\dot \cv_2(t) = 4 \pi r^{n+2} \GG_2(\cc_2(t))$. Now substituting the volume derivatives into \cref{eqn:GROW3} with the identity $\cv_1(t)+\cv_2(t) = \frac{4}{3}\pi r(t)^3$,  we obtain the following growth kinetics describing the change in radius and volume composition for given a radius, composition, and concentrations:
\begin{align}
\ccG(\cc(t),\x)&=\begin{pmatrix}\tfrac{\GG_1(\cc_1(t))+ \GG_2(\cc_2(t))}{\x_1^{-n}}\\ \tfrac{3\GG_1(\cc_1(t)) - 3\left(\GG_1(\cc_1(t))+\GG_2(\cc_2(t))\right)\x_2}{\x_1^{1-n}}\end{pmatrix} \label{eqn:GROW8}.
\end{align}

Note that $\x_1$ and $\x_2$ denote the radius and the volume fraction coordinates, respectively. The characteristics ODE \cref{eqn:CHAR1} with the growth kinetics \cref{eqn:GROW8} can be solved analytically for given time-dependent concentrations $\cc$:

\begin{align}
&\xxxico [t,\x](\tau) \notag\\
&=
\bigg(\!\x_1^{1-n}+(1-n)\!\int_t^{\tau}\!\!\!\!\GG_1(\cc_1(s))+\GG_2(\cc_2(s))\dd s\bigg)^{\!\!\frac{1}{1-n}} \label{eqn:XI1},\\
&\xxxict[t,\x](\tau)\notag\\
&=
3\int_t^{\tau}\tfrac{\GG_1(\cc_1(\omega))}{\xxxico[t,\x](\omega)^{1-n}}e^{-3\int^{\tau}_{\omega}\frac{\left(\GG_1(\cc_1(\eta))+\GG_2(\cc_2(\eta))\right)}{\xxxico[t,\x](\eta)^{1-n}}d\eta}d\omega \nonumber \\
&\qquad +\x_2e^{-3\int_t^{\tau}\frac{\left(\GG_1(\cc_1(\eta))+\GG_2(\cc_2(\eta))\right)}{\xxxico[t,\x](\eta)^{1-n}}d\eta}\label{eqn:XI2},
\end{align}
for $t,\tau \in [0,T]$, and $\x \in \mX:= \R_{\geq \x_n} \times [0,1]$ with minimal particle radius $\x_n \in \R_{>0}$. The fixed-point problem \cref{eqn:COM2} with the analytical solution of the characteristics \cref{eqn:XI1} and \eqref{eqn:XI2} will be discretized in the following section.

\section{Numerical approximation of the fixed point problem}
We can discretize the integro-differential equation by decomposing the time domain: $0 = \ct_1 < \ct_2 < \cdots < \ct_{N_t} = T$ into ${N_t}+1$ time steps, and quadrature points $(\x_\ell)_{\ell=1\ldots,N_t} \subset \mathcal X^N_t$, with quadrature weights $(\oomega_\ell)_{\ell=1\ldots,N_t} \subset \R^{N_t}_{>0}$, such that:
\begin{align}
    \iint_{\mathcal X} q_0(x) V(x) \dd x \approx \sum_{\ell=1}^{N_t} q_0(\x_\ell) \oomega_\ell.
\end{align}
Taking advantage of the semi-group property of the characteristics, i.e., $\xxxic[s,\xxxic[t,\x](s)](\tau) = \xxxic[t,\x](\tau)$ (see e.g. \cite{keimer2017existence}) which also holds for its discretized version $\hxxiC$, we can approach this very efficiently as described in Algorithm \ref{algo:1} using the abbreviations ${\txxxic}(\ell,k):= {\hxxiC}(1,\x_\ell,k) \approx \xxxic (\ct_1,\x_\ell,\ct_k)$ and $\bdelta_k:= \ct_{k+1}-\ct_k$. In Algorithm \ref{algo:1}, the computed $\sC_{i,k}$ is an approximation of the concentration of the $i$-th component at time $t_k$, i.e. $\sC_{i,k} \approx \cc_i(\ct_k)$.
\RestyleAlgo{ruled}
\begin{algorithm*}[h]
\label{algo:1}
\caption{Numerical approximation scheme for the integral-fixed-point equation \cref{eqn:COM2} which prescribes the evolution of the driving forces over process time.}
~\\[5pt]
\textbf{initialization step}\\[20pt]

\For{$\ell = 1 \ldots Nx$}{
$\txxxic (\ell,1) \gets \x_{\ell} $ 
}
$\sC_{1,1}  \gets \tfrac{1}{\mathcal V} \m_1(\ct_{1})-\tfrac{\brho_1}{\mathcal V}\sum_{\ell=1}^{N_x}\V_1\left(\x_\ell\right)q_0(\x_\ell)\oomega_\ell $

$\sC_{2,1}  \gets \tfrac{1}{\mathcal V} \m_2(\ct_{1})-\tfrac{\brho_2}{\mathcal V}\sum_{\ell=1}^{N_x}\V_2\left(\x_\ell\right)q_0(\x_\ell)\oomega_\ell $\\[20pt]

\textbf{time-stepping}\\[10pt]

\For{$k = 1 \ldots N_t-1$}{
\For{$\ell = 1 \ldots Nx$}{
$\txxxic (\ell,k+1) \nonumber \gets  \begin{pmatrix}
\left(\txxxico(\ell,k)^{1-n} + (1-n)\bdelta_k\left(\GG_1(\sC_{1,k}) + \GG_2(\sC_{2,k})\right)\right)^{\frac{1}{1-n}} \\
 \left(3\tfrac{\GG_1(\sC_{1,k})\bdelta_k}{\txxxico(\ell,k)^{1-n}} + \txxxict(\ell,k)\right)\exp\left(-3\tfrac{\left(\GG_1(\sC_{1,k})+\GG_2(\sC_{2,k})\right)\bdelta_k}{\txxxico(\ell,k)^{1-n}} \right) \end{pmatrix}  $
 }
$\sC_{1,k+1} \nonumber \gets \tfrac{1}{\mathcal V} \m_1(\ct_{k+1})-\tfrac{\brho_1}{\mathcal V}\sum_{\ell=1}^{N_x}\V_1\left(\txxxic(\ell,k+1)\right)q_0(\x_\ell)\oomega_\ell $ 

$\sC_{2,k+1} \nonumber \gets \tfrac{1}{\mathcal V} \m_2(\ct_{k+1})-\tfrac{\brho_2}{\mathcal V}\sum_{\ell=1}^{N_x}\V_2\left(\txxxic(\ell,k+1)\right)q_0(\x_\ell)\oomega_\ell $ 
}
\end{algorithm*}

As the complexity of every iteration step is of order $N_x$ and we have $N_t$ steps, this approach yields the order $\mathcal O(N_x N_t)$ in the numerical scheme.

\begin{figure*}[h!]
    \centering
{\includegraphics[scale=0.65,clip=true,trim=0 0 22 0]{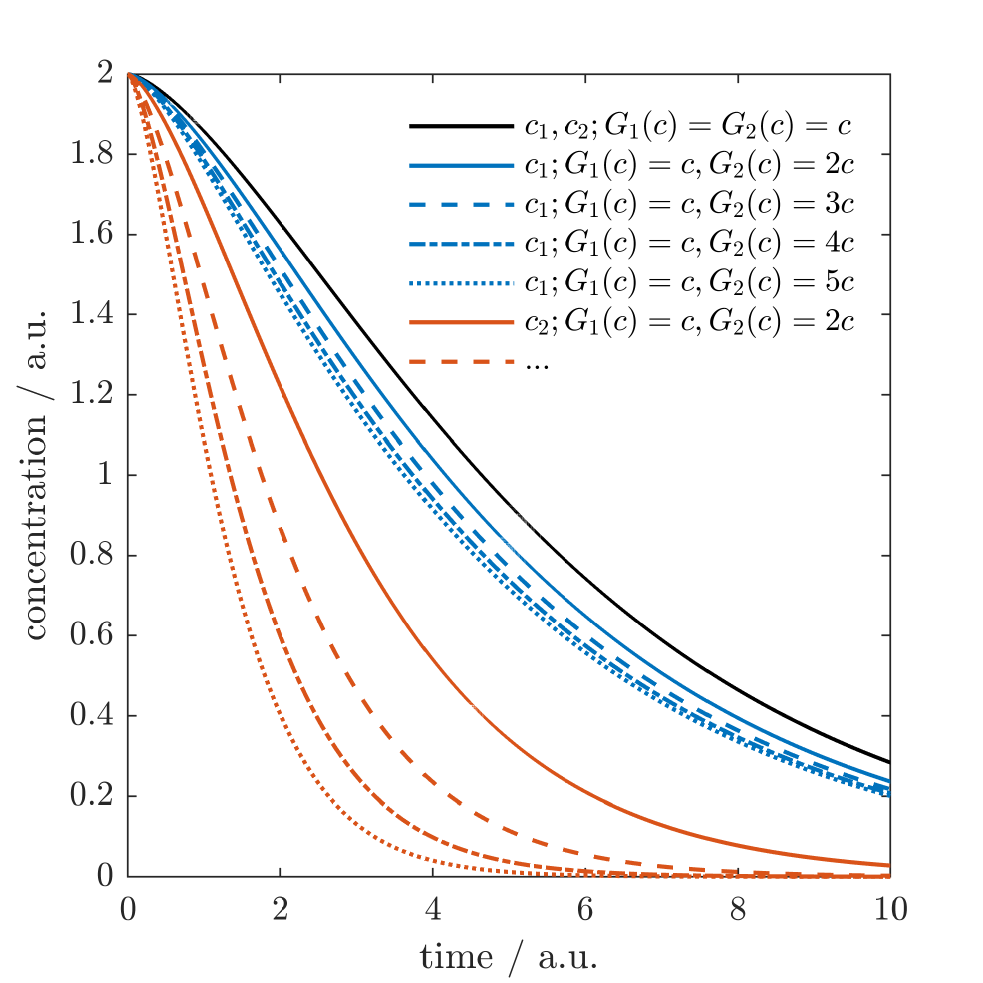}}
{\includegraphics[scale=0.65,clip=true,trim=0 0 25 0]{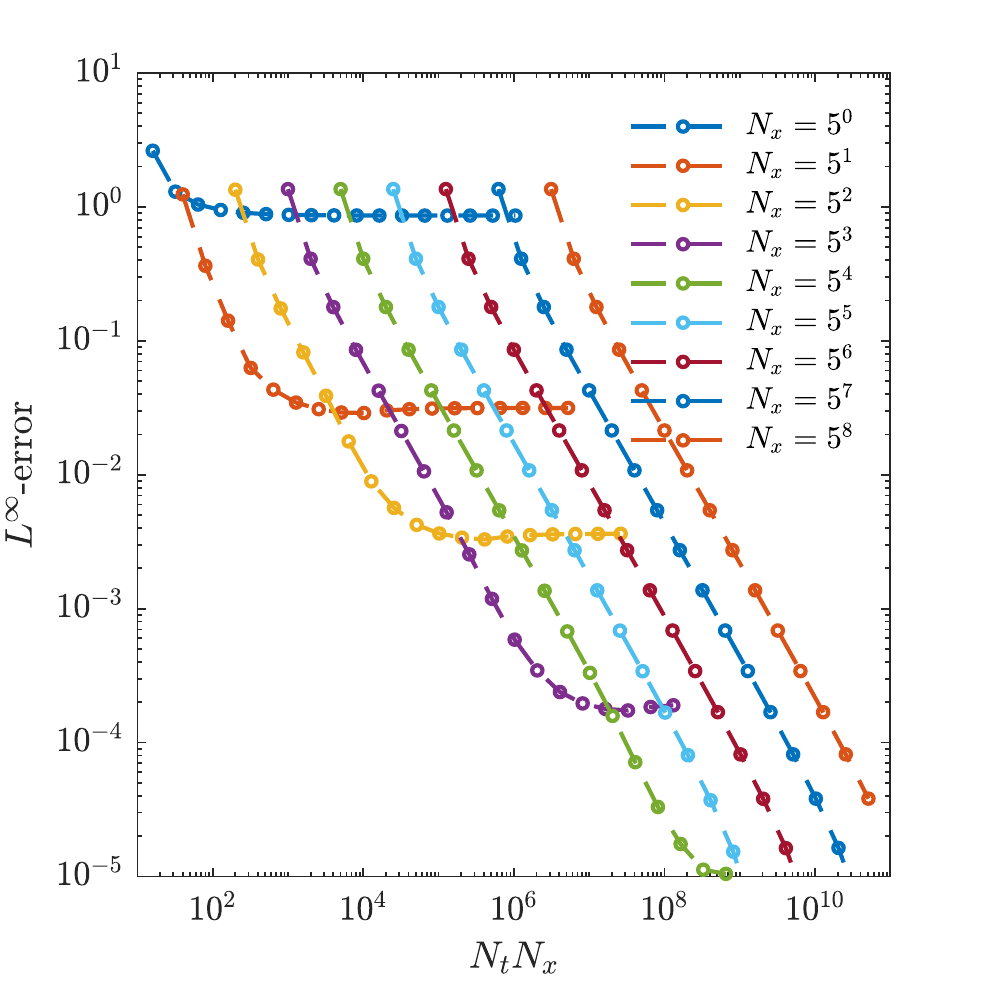}}
{\includegraphics[scale=0.65,clip=true,trim=0 0 25 0]{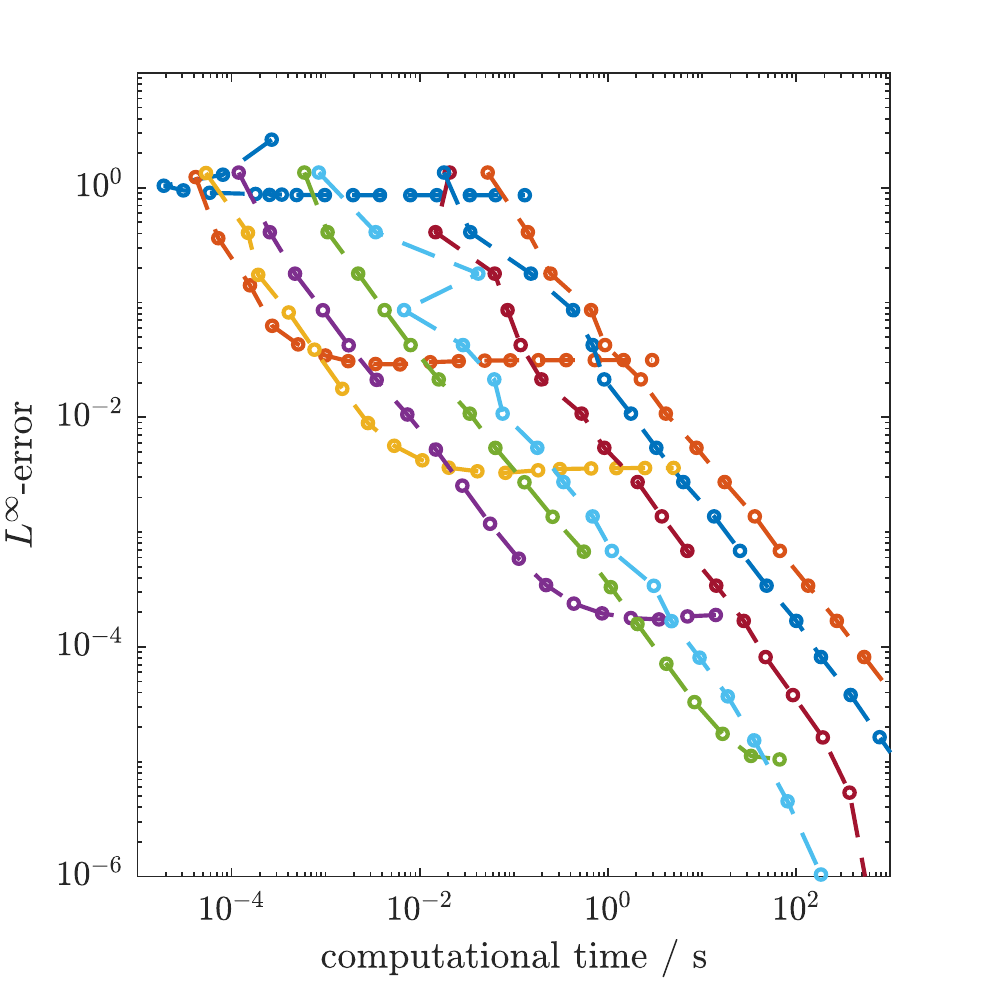}}
\caption{Solute concentrations over process time for different growth kinetics (\textbf{left}); Numerical approximation of the maximal absolute error in the concentrations ($L^\infty$-error) for different time discretizations $N_t$, and disperse property discretizations of the disperse property as a function of the number of degrees of freedom (DoF) $N_tN_x$ (\textbf{middle}); and the total computational time (\textbf{right}). The error analysis is performed for $\GG_1(c) = c, \GG_2(c) = 5c$. The reference solution is given by $N_t \approx 10^6, \& \ N_x \approx 4 \; 10^5$. For fixed $N_x$, we obtain a clear first-order behavior with respect to both the DoF and to the computational time, provided the spatial resolution is large enough.}
\label{fig:convergence}
\end{figure*}

To study the numerical performance of MD-eMoM, we conduct a numerical error analysis by refinement series. We assume the following growth kinetics: $\GG_1(c)=c$, 
$\mathcal V=1$, $\cc_1(0) = \cc_2(0)=2$, 
$\brho_1=\brho_2=1$,
$\V_1(\x)= \frac{4}{3} \pi  \x_1^3  \x_2$, and 
$\V_2(\x)= \frac{4}{3} \pi  \x_1^3  (1-\x_2)$; and take as initial datum
$q_0(\x) = \max\big\{1-(\x_1- \bar \x_1)^2\delta_1^{-2} + (\x_2- \bar \x_2)^2\bdelta_2^{-2} \ , 0 \big\}^2$, with $\bar \x = (0.1,0.75)^\top$ and $\bdelta = (0.05,0.25)^\top$. To emphasize the impact of different growth kinetics on the composition evolution, we consider different growth rates for the second educt component ($\GG_2$). The effect of different growth rate ratios $\GG_2/\GG_1$ on the temporal concentration evolution of each component is shown in \Cref{fig:convergence} (leftmost figure). The numerical approximation error w.r.t. changes in the discretization is shown in the middle (on the basis of the degrees of freedom $N_xN_t$) and right plot (on the basis of the computational time). The computational time refers to a vectorized MATLAB R2021b implementation of Algorithm \ref{algo:1} on a MacBook Pro 2021, Apple M1 Chip with 32GB RAM. In the same manner, higher-order explicit or implicit discretization schemes can easily be derived.


\section{Numerical approximation of the population balance equation solution}
The PBE solution can be accessed in two different ways based on the solution of the fixed-point equation. Firstly, the particle density function for a given disperse property $\x \in \mX$ at a given time $t$, i.e., $q(t,\x)$, is already given by \cref{eqn:SOL1}. Secondly, the solution can be evaluated with the introduced characteristics $\xxxic$, i.e.:
\begin{align}
q(t,\xxxic[0,\x](t))=q_0(\x)\ \text{det}\left(D_2\xxxic[t,\xxxic[0,\x](t)](0)\right)\label{eqn:SOL3}.
\end{align}
To evaluate both solution formulas, we approximate the characteristics stated in \cref{eqn:XI1} and \eqref{eqn:XI2} by discretizing the integrals into the time interval and then approximating the integrands by the interval's value at the lower time bound, i.e.: 

\begin{equation}
\begin{split}
&\xxxico[\ct_k,\x](\ct_i)\\
&\quad\approx 
\Big(\x_1^{1-n}+(1-n)\textstyle\sum_{\ell=k}^{i-1} \textbf{G}_\ell\Big)^{\frac{1}{1-n}} \\
&\quad=: \hxxiCo(k,\x,i),\\[10pt]
&\xxxict[\ct_k,\x](\ct_i)\\
&\quad\approx 3\sum_{\ell=k}^{i-1} \frac{\GG_1(\sC_{1,\ell})(\ct_{\ell+1}-\ct_\ell) e^{-3\sum_{j=\ell}^{i-1}\frac{ \textbf{G}_j}{\hxxiCo(k,\x,j)^{1-n}}}}{\hxxiCo\!(k,\x,\ell)^{1-n}}\\
&\qquad\qquad+\x_2 e^{-3\sum_{\ell=k}^{i-1}\frac{ \textbf{G}_\ell}{\hxxiCo(k,\x,\ell)^{1-n}}}\\
&\quad=: \hxxiCt(k,\x,i),
\end{split}
\end{equation}
with 
\begin{align*}
 \textbf{G}_\ell := \big(\GG_1(\sC_{1,\ell})+\GG_2(\sC_{2,\ell})\big)\big(\ct_{\ell+1}-{\ct_\ell}\big)
\end{align*}
The above equation holds for $k \leq i$. If $k > i$, the summations have to be understood as follows: $\sum_{\ell=k}^{i-1}=-\sum_{\ell=i}^{k-1}$. The solution of the PBE can then be approximated as follows:

\begin{equation}
\begin{split}
q(\ct_k,\x) \approx q_0\left(\hxxiC(k,\x,1)\right)\left(\frac{
\hxxiCo(k,\x,1)}{\x_1}\right)^n  \Psi_{k,k}(x), \\
q\left(\ct_k,\hxxiC(1,\x,k)\right) \approx q_0(\x) \left(\frac{\x_1}{
\hxxiCo(1,\x,k)}\right)^n \Psi_{1,k}(x).
\end{split}
\end{equation}
with 
\begin{align*}
     \Psi_{i,k}(x) := \exp\left(3\sum_{\ell=1}^{k-1}\frac{\left(\GG_1(\sC_{1,\ell})+\GG_2(\sC_{2,\ell})\right)\left(\ct_{\ell+1}-\ct_\ell\right)}{\hxxiCo(i,\x,\ell)^{1-n}} \right)
\end{align*}

\begin{figure*}[h!]
\centering
{\includegraphics[scale=0.70,clip,trim=70 0 70 0]{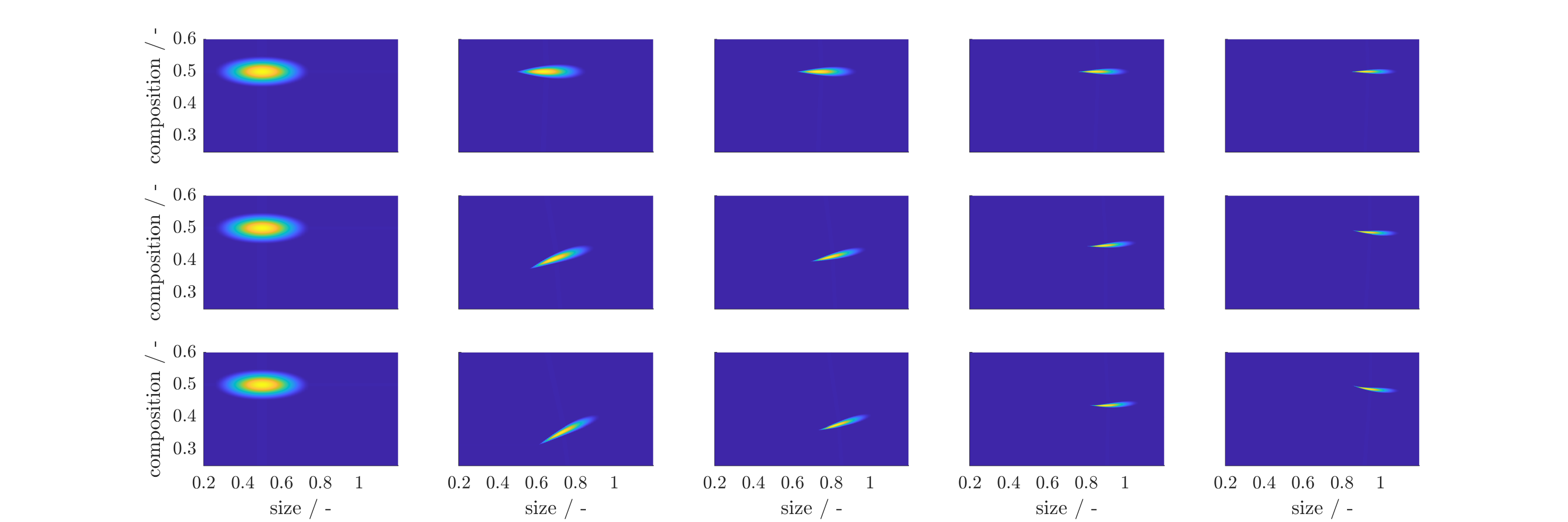}}
\caption{Evolution of the PSD for $G_1\equiv G_2$ (\textbf{top row}), $G_1\equiv 2G_2$ (\textbf{middle row})  and $G_1 \equiv 5 G_2$ (\textbf{bottom row}) from left to right. The predicted PSD focusing, denoted by the growth kinetics, clearly shows the strong ability of eMoM to capture this.}
\label{fig:PSD_evolution}
\end{figure*}

The first solution formula allows us to evaluate the particle distribution at any time, and disperse property, while the second enables the initial datum to be easily tracked over time, as only the disperse properties $\x$ within the support of the initial datum $q_0$ have to be evaluated. It is worth mentioning that whatever scheme is used to approximate the evolution of the concentrations over time, the solution formula can then be utilized to evaluate the PBE solution with high accuracy -- keeping in mind that this is only true if the concentrations are approximated accurately enough. In \Cref{fig:PSD_evolution}, the evolution of a PSD is depicted for two growth rate ratios, showing a clear change in the composition over process time.



\subsection{Comparison with classical discretization schemes}
\begin{figure*}
    \centering
\begin{minipage}[b]{0.36\linewidth}
{\includegraphics[scale=0.70,clip=true,trim=0 0 0 0]{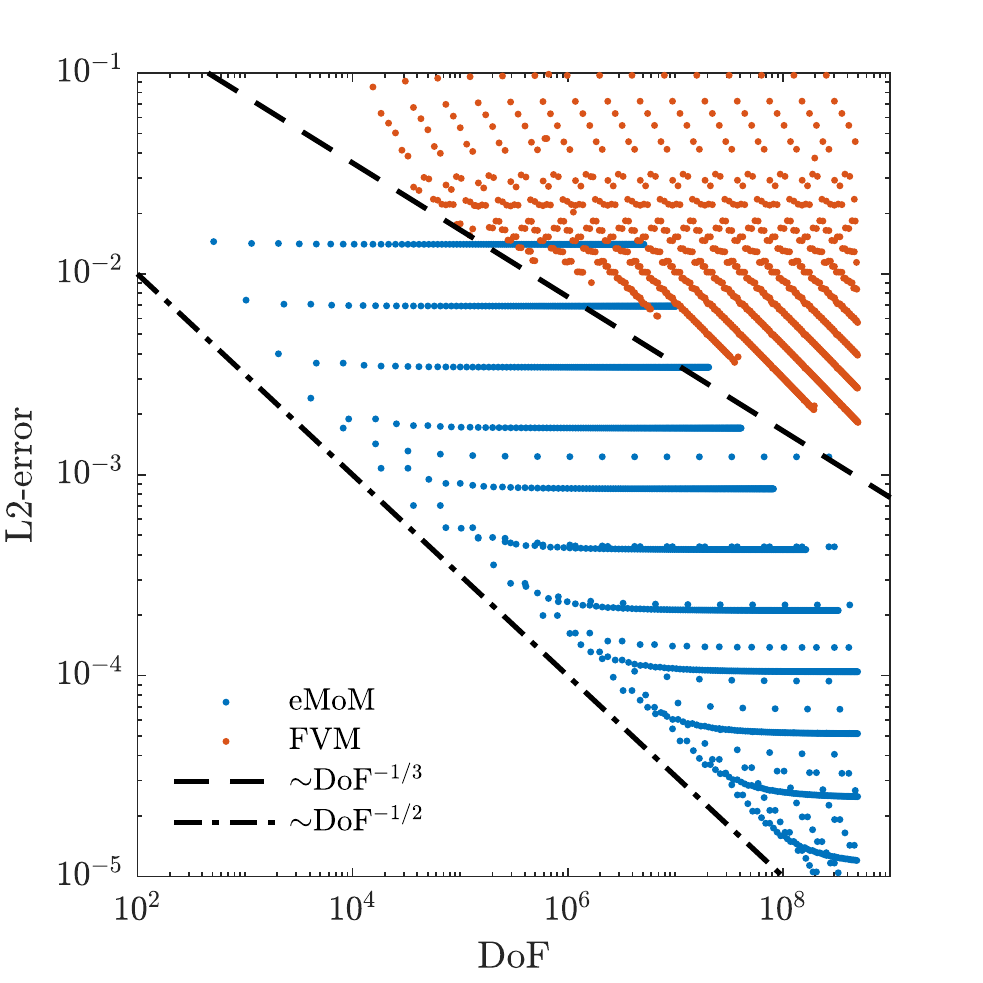}}
\end{minipage}
\begin{minipage}[b]{0.63\linewidth}
{\includegraphics[scale=0.70,clip=true,trim=25 40 15 10]{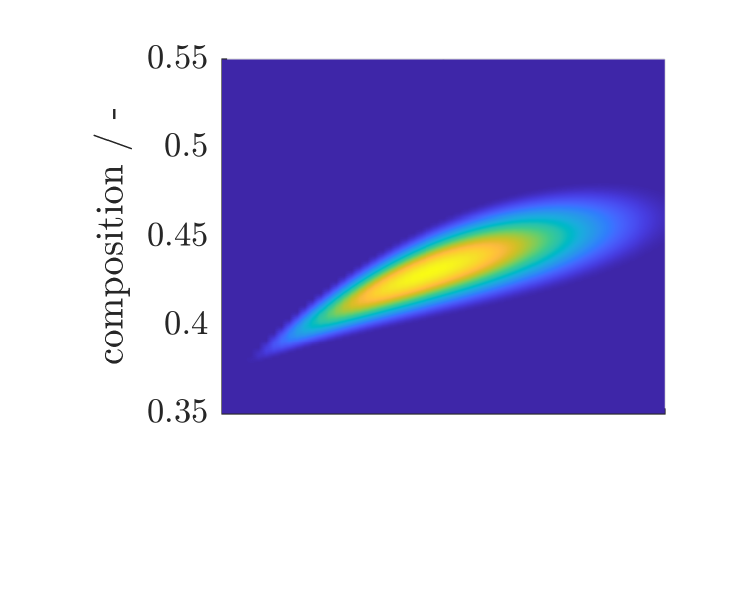}}
{\includegraphics[scale=0.70,clip=true,trim=57 40 15 10]{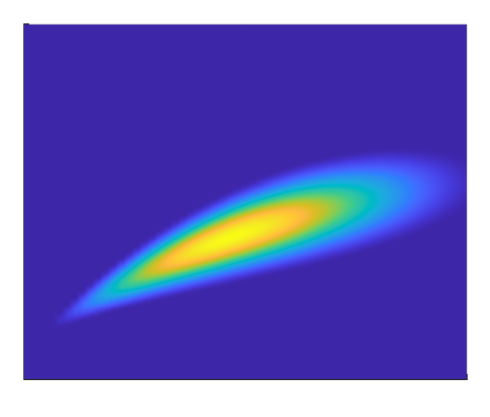}}
\includegraphics[scale=0.68,clip=true,trim=
57 40 15 10]{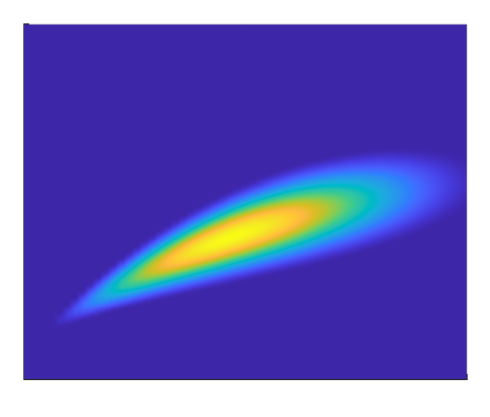}

{\includegraphics[scale=0.70,clip=true,trim=25 20 15 2]{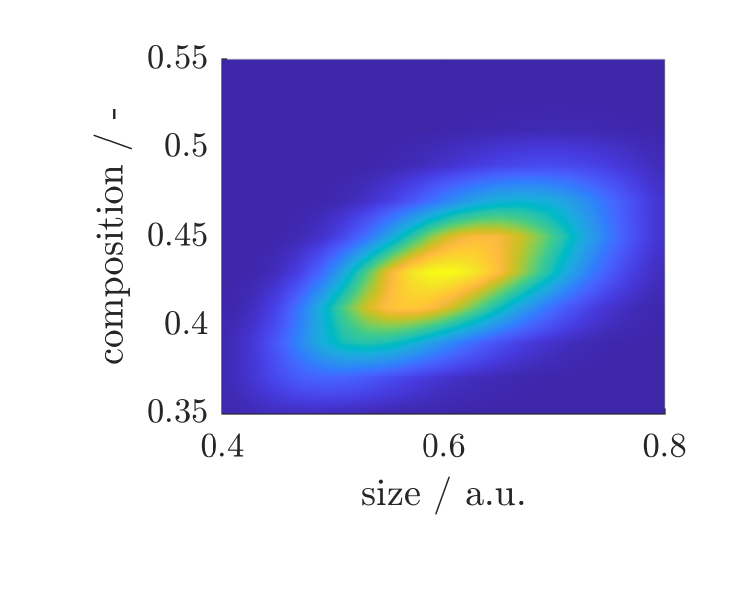}}
{\includegraphics[scale=0.70,clip=true,trim=57 20 15 2]{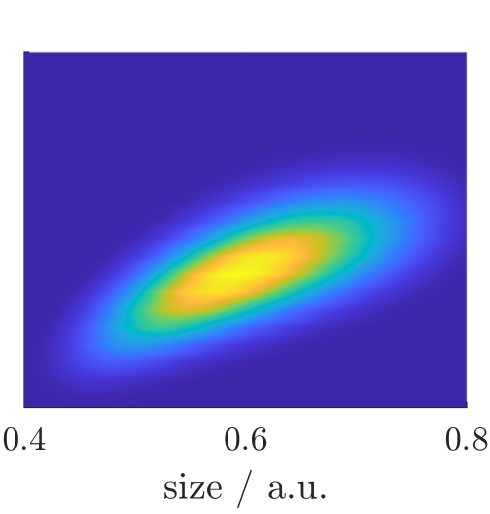}}
{\includegraphics[scale=0.70,clip=true,trim=57 20 15 2]{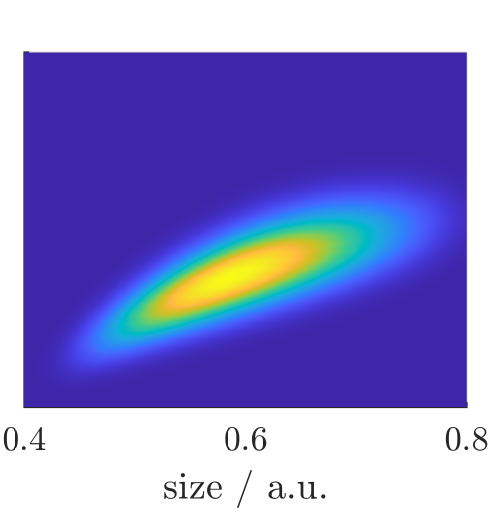}}
\end{minipage}

\caption{\textbf{Left:} Comparison of the $L^2$-error for the predicted evolution of the concentration between MD-eMoM and FVM for a large range of discretizations (DoF). One can clearly identify an order of $1/3$ with respect to the DoF for FVM and of $1/2$ for eMoM. \textbf{Right:} The number density function $q$ at $t=0.01$ for different DoF (from left to right, $10^6$ through $10^7$ to $10^8$), where the top and bottom row represent the solution of MD-eMoM and FVM, respectively. For coarse discretizations (small DoF), the FVM solution is diffusive. Furthermore, the whole domain in FVM needs to discretized though it is only partly used. By contrast, eMoM is non-diffusive and, by definition, only implicitly discretizes the support of the PSD for all process times.}
    \label{fig:eMoM_fvm}
\end{figure*}

We compare the numerical scheme derived here with state-of-the-art finite volume method schemes (FVM) \cite{leveque2002finite,braatz2002particle}. To enforce boundedness and numerical stability of the FVM, we consider the total variation diminishing schemes (TVD) \cite{leveque2002finite} with a \textit{van-Leer} limiter function enabling weighting between first and second order approximations of the growth term\cite{Sweby_TVD}. \Cref{fig:eMoM_fvm} shows a numerical comparison between MD-eMoM and FVM. The newly derived MD-eMoM clearly outperforms FVM. First, it does not need to satisfy a CFL condition \cite{courant1928partiellen} because it relies on the solution of the characteristics. Second, MD-eMoM does not suffer from poor (coarse) discretization of the disperse property as only the support of the initial datum has to be discretized, see \Cref{fig:eMoM_fvm}.

\section{Additional insights into Nanoparticle structure}
\begin{figure*}[h!]
\centering
{\includegraphics[scale=0.7,clip,trim=0 0 0 10]{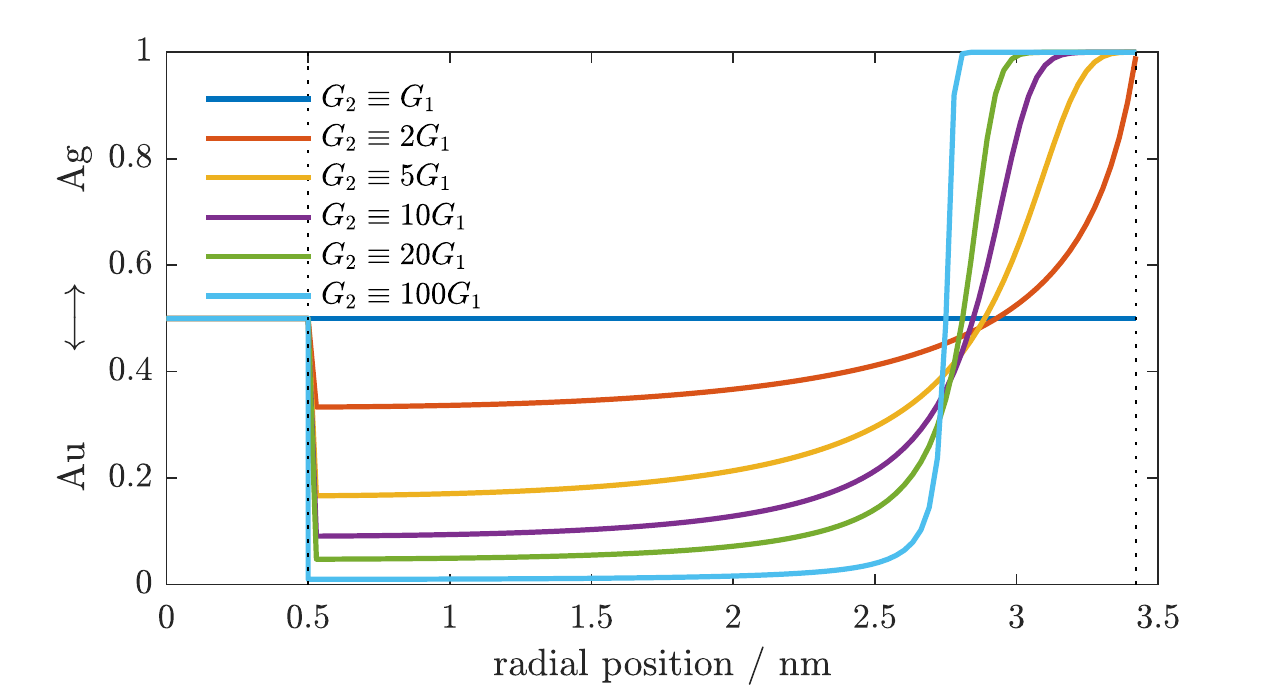}}
\includegraphics[scale=0.7,clip,trim=0 0 0 10]{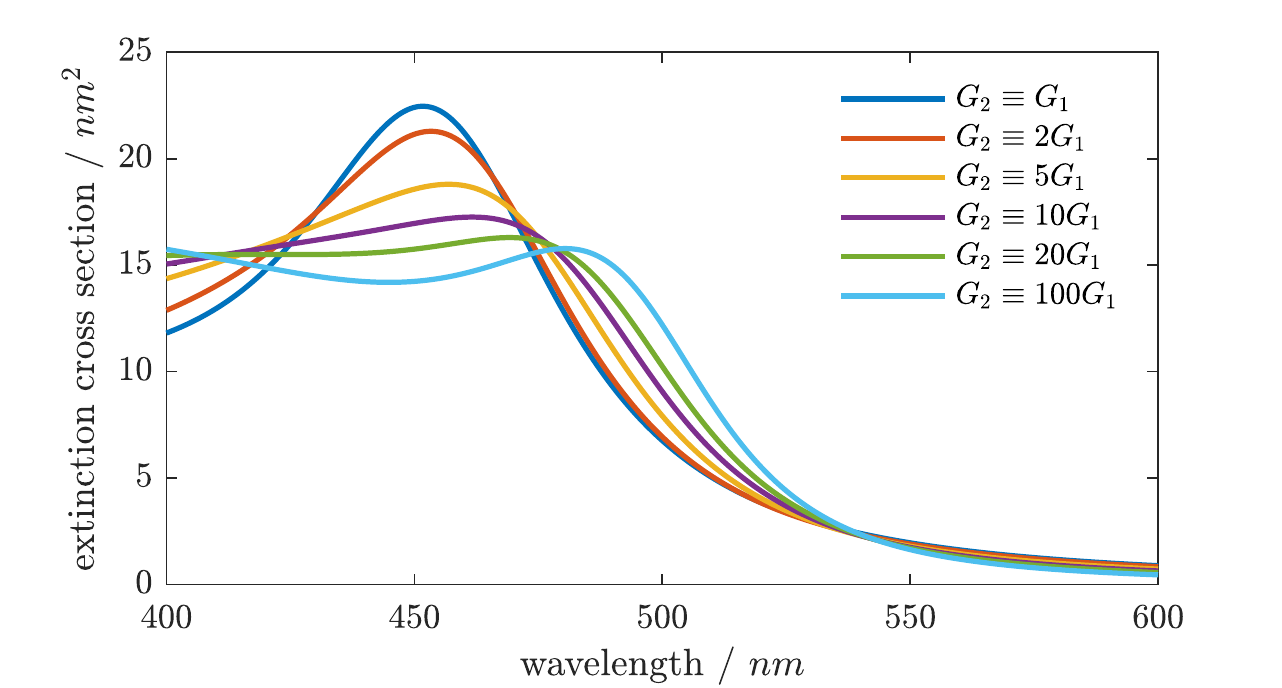}
\caption{Radial composition (\textbf{left}) and the resulting extinction cross-section spectra (\textbf{right}) of AuAg nanoalloys for different growth rate ratios. The optical properties heavily depend on the radial distribution of Au-Ag. Note that all the particles have the same average fraction of Au-Ag, i.e., 0.5. For the sake of simplicity, the initial particle size distribution $q_0$, i.e., the seeds, was chosen as a Dirac-$\delta$ with a radius of $0.5$nm and a composition of $0.5$. For the growth rate ratios $1$ and $5$, the evolution of a full $q_0$ is depicted in \Cref{fig:PSD_evolution}. The Dirac follows the same evolution.}
    \label{fig:radial_comp}
\end{figure*}

In many applications, the inner-particle structure determines the product properties, e.g., the radial composition of gold (Au) and silver (Ag) in AuAg nanoalloys determine the optical properties \cite{rioux2015seeded}. Having the time evolution of the educt concentrations $\cc$, we can also reconstruct the inner-particle composition for every particle in the number density function based on the characteristics. For a particle at time $t$ with initial disperse property $\x^0$, the composition $\mathcal F$ at every radial position $\xxxic(0,\x^0,\tau)$ for $\tau \in [0,t]$ is given by:
\begin{align}
\mathcal F_{\x^0,t}\big(\xxxico(0,\x^0,\tau)\big) =  \frac{\GG_1(\cc_1(\tau))}{\GG_1(\cc_1(\tau))+\GG_2(\cc_2(\tau))}.
\end{align}
This allows us to trace the evolution of the radial composition over process time and to characterize the properties of the final particle size distribution. \\
To demonstrate the potency of this unique analysis approach, we consider the seeded growth of gold-silver alloy NPs and investigate the effect of different growth rate ratios on the inner-particle structure (here optical properties). Under the assumption of radial symmetry, the optical properties can be numerically calculated by MIE Theory \cite{mie1908beitrage} (we use the MATLAB code of J. Schäfer\cite{schafer2011implementierung}). As composition-dependent material properties, i.e., refractive indices, we interpolate the measured refractive index for gold-silver alloys by McPeak \textit{et. al.} \cite{mcpeak2015plasmonic}. We prescribe the initial datum $q_0$ to be of uniform composition 0.5. \Cref{fig:radial_comp} (left) shows the inner-particle composition for a particle with size 7 nm after the solid formation process is finished. The different radial compositions due to the different growth rate ratios result in different extinction spectra, and thus optical properties, as seen in \Cref{fig:radial_comp} (right).

\section{Conclusion and Outlook}
The multi-dimensional exact method of moment is an efficient way to approximate solutions of multi-dimensional population balance equations by relying on a fixed-point reformulation in terms of the driving forces -- here the solute concentrations. The introduced numerical approach allows for highly accurate prediction of the evolution of multi-dimensional number density function, and the reformulation of the governing equations into a set of integro-differential equations results in markedly improved computational efficiency and numerical accuracy compared to state-of-the-art finite volume or finite element methods. Another advantage relates to process optimization, as the derived numerical scheme is differentiable by construction and thus derivatives with respect to process conditions can be computed by utilizing the implicit-function theorem.

Our numerical scheme can easily be extended to take into account the nucleation and growth of composite particles, as well as the formation of anisotropic particles\cite{braatz2002particle}. Moreover, it is not limited to just one composition but can handle multiple compositions. A natural straightforward extension of the scheme would be implicit discretization of the fixed-point problem, resulting in an $n$-dimensional nonlinear system of equations for each time-step in Algorithm \ref{algo:1}, ($n$ being the number of considered concentrations). Furthermore, the multi-dimensional exact method of moment idea can be coupled to fluid-flow, as already outlined in Bänsch et. al.\cite{bansch2022highly} for the one-dimensional case.

\section*{Acknowledgements}
L. Pflug has been supported by the Deutsche Forschungsgemeinschaft (DFG, German Research Founda- tion) – Project-ID 416229255 – SFB 1411.

\bibliographystyle{model1-num-names}
\bibliography{biblioNEW2.bib}

\end{document}